\documentclass{amsart}

\newtheorem{theorem}{Theorem}
\newtheorem{lemma}{Lemma}

\newcommand{\bt}{\begin{theorem}}
\newcommand{\et}{\end{theorem}}
\newcommand{\bl}{\begin{lemma}}
\newcommand{\el}{\end{lemma}}
\newcommand{\beq}{\begin{equation}}
\newcommand{\eeq}{\end{equation}}

\usepackage{amsmath,amssymb}
\usepackage[all]{xy}

\title{Infinite sumsets with many representations}
\author{Melvyn B. Nathanson}
\address{Department of Mathematics\\Lehman College (CUNY)\\Bronx, NY 10468}

\email{melvyn.nathanson@lehman.cuny.edu}

\dedicatory{To Helmut Maier on his 60th birthday}

\thanks{This work was supported in part by grants 
from the PSC-CUNY Research Award Program.}

\date{\today}

\begin{document}

\maketitle

\begin{abstract}
Let $A$ be an infinite set of nonnegative integers.  For $h \geq 2$, 
let $hA$ be the set of all sums of $h$ not necessarily distinct elements of $A$.  
If every sufficiently large integer in the sumset $hA$ has at least two 
representations, then $A(x) \geq (\log x)/\log h)-w_0$, 
where $A(x)$ counts the number of integers $a \in A$ such that $1 \leq a \leq x$.
\end{abstract}

\section{Representation functions of sumsets}
Let $A$ be a set of integers and let $h \geq 2$ be an integer.  
The \emph{counting function} $A(x)$  counts the number of positive 
integers in the set $A$ that do not exceed $x$.  
The $h$-fold sumset $hA$ is the set of all integers $n$ 
that can be written as sums of $h$ not necessarily distinct elements of $A$.  
For every  integer $n$,  the 
\emph{representation function}
$r_{A,h}(n)$ counts the number of $h$-tuples 
$(a_1, a_2, \ldots, a_h) \in A^h$ 
such that 
\[
a_1 \leq a_2 \leq \cdots \leq a_h
\] 
 and 
 \[
 a_1 +a_2 + \cdots + a_h = n.
 \]
 
A \emph{Sidon set} is a set $A$ of nonnegative integers such that every element 
in the sumset $2A$ has a unique representation, that is, 
$r_{A,2}(n) = 1$ for all $n \in 2A$.
More generally, for positive integers $h$ and $s$, a $B_{h,s}$-set is 
a set $A$ of  nonnegative integers such that 
$r_{A,h}(n) \leq s$ for all $n \in hA$.
Sets whose sumsets have few representations 
have been studied intensively (cf. Halbertam-Roth~\cite{halb-roth66}, 
O'Bryant~\cite{obry04}).

In this paper we consider sets whose $h$-fold sumsets have many representations.
A basic result is that if $h \geq 2$ and $A$ is an infinite set of nonnegative integers 
with $r_{A,h}(n) \geq 2$ for all sufficiently large 
integers $n \in hA$, then 
\[
A(x) \gg \log x.  
\]
In the special case $h=2$, Balasubramanian and Prakesh~\cite{bala-prak04} 
proved that there is a number $c > 0$ such that, if $A$ is an infinite set of 
nonnegative integers with $r_{A,2}(n) \geq 2$ for all sufficiently large 
integers $n \in 2A$, then 
\[
A(x) \geq c\left(\frac{\log x}{\log\log x} \right)^2.
\]
This improved a previous result of 
Nicolas, Ruzsa, and S\" ark\` ozy~\cite{nico-ruzs-sark98}, 
who also proved the existence of an infinite set $A$ of 
nonnegative integers with $r_{A,2}(n) \geq 2$ for all sufficiently large 
integers $n \in 2A$ such that 
\[
A(x) \ll (\log x)^2.
\] 
It is an open problem to extend these results to $h$-fold sumsets  for $h \geq 3$.  

\emph{Acknowledgements.}  I thank Michael Filaseta for  bringing these problems to my attention, and  Quan-Hui Yang 
for the reference to the paper of Balasubramanian and Prakesh.

\section{Growth of sets with many representations}

Let $[u,v)$ denote the interval  of integers $i$ such that $u \leq i < v$.
Let $|X|$ denote the cardinality of the set $X$.

\bt         \label{filaseta:theorem:h-ell2}
Let $h \geq 2$  be an integer, 
and let $A$ be an infinite set of nonnegative integers.
If $r_{A,h}(n) \geq 2$ for all sufficiently large integers $n \in hA$, 
then there is a positive number $w_0$ such that 
\[
A(x) > \frac{\log x}{\log h} - w_0
\]
for all $x \geq h$.  
\et

\begin{proof}
For every positive integer $k$, let 
\[
I_k = \left[ h^{k-1}, h^k \right)
\]
and
\[
A_k = A \cap I_k.
\]
The sets $\{A_k:k=1,2,\ldots \}$ partition  $A\setminus \{ 0\}$.

There exists a positive integer $n_0$ such that, if $n \geq n_0$ and $n \in hA$, 
then $r_{A,h}(n) \geq 2$.  Because $A$ is infinite, there exists $a_0 \in A$ 
with $ha_0 \geq n_0$.  Choose $k_0$ such that $a_0 \in A_{k_0}$.

Suppose that $k \geq k_0$ and $A_k \neq \emptyset$.  
Let 
\[
a_k^* = \max(A_k).
\]
Then 
\[
h^{k-1} \leq a_k^* < h^k
\]
and 
\[
A \cap \left[a_k^*+1,h^k \right) = \emptyset.
\]
Consider the integer
\[
ha_k^* \in hA.
\]
Because $a_k^* \geq a_0$, we have  $ha_k^* \geq ha_0 \geq n_0$, 
and so $r_{A,h}(ha_k^*) \geq 2$.   
It follows that the set $A$ contains nonnegative integers $a_1,\ldots, a_h$  such that 
\beq   \label{filaseta:ineq1} 
a_1 < a_h
\eeq
\beq   \label{filaseta:ineq2} 
a_1 \leq a_2 \leq \cdots \leq a_h
\eeq
and
\[
 a_1 +a_2 + \cdots + a_h = ha_k^*.
 \]
Because $A$ is a set of nonnegative integers, we have  
\[
a_h \leq ha_k^*.
\]
Inequalities~\eqref{filaseta:ineq1} and~\eqref{filaseta:ineq2}  imply that 
\[
ha_k^*=  a_1 +a_2 + \cdots + a_h < ha_h
\]
and so
\[
a_k^* < a_h  \leq ha_k^* < h^{k+1}.
\]
Therefore,
\[
a_k^* + 1 \leq a_h < h^{k+1}.
\]
Equivalently,  
\[
a_h \in \left[ a_k^*+1,h^{k+1} \right) =  \left[ a_k^*+1,h^k \right) \cup A_{k+1}.
\]
Because $A \cap \left[a_k^*+1,h^k \right) = \emptyset$,
we see that $a_h \in A_{k+1}$
and so $A_{k+1} \neq \emptyset$.  
It follows by induction that $A_k \neq \emptyset$ for all $k \geq k_0$.

Let $x \geq h$, and choose the positive integer $t$ such that 
\[
h^t\leq x < h^{t+ 1}.
\]
Because
\[
A \setminus \{ 0\} = \bigcup_{k= 1}^{\infty} A_k
\]
it follows that 
\[
A(x) \geq A\left( h^t \right) \geq t - (k_0 -1) > \frac{\log x}{\log h} - k_0
\]
for all $x \geq h$.
Let $w_0 = k_0$.  This completes the proof.  
\end{proof}

\bt
Let $\ell \geq 2$  be an integer, 
and let $A$ be an infinite set of nonnegative integers.
If $r_{A,2}(n) \geq \ell$ for all sufficiently large integers $n \in 2A$, 
then there is a positive number $w_0$ such that 
\[
A(x) > \frac{(\ell - 1)\log x}{\log 2} - w_0
\]
for all $x \geq 2$.
\et

\begin{proof}
For every positive integer $k$, let 
\[
I_k = \left[ 2^{k-1}, 2^k \right)
\]
and
\[
A_k = A \cap I_k.
\]
The sets $\{A_k:k=1,2,\ldots \}$ partition  $A\setminus \{ 0\}$.

There exists a positive integer $n_0$ such that if $n \geq n_0$ and $n \in 2A$, 
then $r_{A,2}(n) \geq \ell$.  Because $A$ is infinite, there exists $a_0 \in A$ 
with $2a_0 \geq n_0$.  Choose $k_0$ such that $a_0 \in A_{k_0}$.

Suppose that $k \geq k_0$ and $A_k \neq \emptyset$.  
Let 
\[
a_k^* = \max(A_k).
\]
Then 
\[
2^{k-1} \leq a_k^* < 2^k
\]
and 
\[
A \cap \left[a_k^*+1,2^k \right) = \emptyset.
\]
Consider the integer
\[
2a_k^* \in 2A.
\]
Because $a_k^* \geq a_0$, we have  $2a_k^* \geq 2a_0 \geq n_0$, 
and so $r_{A,2}(2a_k^*) \geq \ell$.   
It follows that the set $A$ contains a subset 
\[
\{ a_{i,j}:i=1,2 \text{ and } j = 1,\ldots, \ell -1\}
\]
such that,  for $j=1,2,\ldots, \ell -1$,  
\[
a_{1,j} < a_{2,j} 
\] 
\[
a_{1,1} <  a_{2,1} < \cdots < a_{\ell -1,1} < a_k^* 
<  a_{\ell-1,2}< \cdots < a_{2,2} < a_{1,2} \]  
and
\[
 a_{1,j} +a_{2,j} = 2a_k^*
 \]
 for $j=1,2,\ldots, \ell-1$.   
Moreover,  $a_{1,j}  \geq 0$ implies that 
\[
a_{2,j} \leq 2a_k^*
\]
for $j =1,2,\ldots, \ell-1$.    
 
Because $0 \leq a_{1,j} < a_{2,j} $, we have 
\[
2a_k^*=  a_{1,j} + a_{2,j} < 2a_{i,2}
\]
and so
\[
a_k^* < a_{2,j}  \leq 2a_k^* < 2^{k+1}.
\]
Therefore,
\[
a_k^* + 1 \leq a_{2,j} < 2^{k+1}.
\]
Equivalently, 
\[
a_{2,j}  \in \left[ a_k^*+1,2^{k+1} \right) =  \left[ a_k^*+1, 2^k \right) \cup A_{k+1}
\]
 for $j=1,2,\ldots, \ell-1$.  
Because $A \cap \left[a_k^*+1, 2^k \right) = \emptyset$,
we see that $a_{2,j} \in A_{k+1}$ 
 for $j =1,2,\ldots, \ell-1$,
and so $\left| A_{k+1}\right| \geq \ell - 1$.  
It follows by induction that  $\left| A_k \right| \geq \ell - 1$ for all $k \geq k_0 + 1$.

Let $x \geq 2$, and choose the positive integer $t$ such that 
\[
2^t \leq x < 2^{t + 1}.
\]
Because
\[
A \setminus \{ 0\} = \bigcup_{k= 1}^{\infty} A_k
\]
it follows that 
\begin{align*}
A(x) & \geq A\left( 2^t \right) \geq  (\ell - 1)(t - k_0 ) \\
& > (\ell - 1) \left( \frac{\log x}{\log 2} - k_0 - 1 \right) \\
& = \frac{(\ell - 1)\log x}{\log 2} - w_0
\end{align*}
for all $x \geq 2$.  
Let  $w_0 =  (\ell - 1) (k_0 + 1)$.  
This completes the proof.  
\end{proof}

\bt
Let $h \geq 2$ and $\ell \geq 2$ be integers, 
and let $A$ be an infinite $B_{h-1, s}$ set of nonnegative integers.   
If $r_{A,2}(n) \geq \ell$ for all sufficiently large integers $n \in hA$, 
then there is a positive number $w_0$ such that 
\[
A(x) >   \frac{(\ell - 1)\log x}{s \log h}  - w_0
\]
for all $x \geq h$.
\et

\begin{proof}
For every positive integer $k$, let 
\[
I_k = \left[ h^{k-1}, h^k \right)
\]
and
\[
A_k = A \cap I_k.
\]
The sets $\{A_k:k=1,2,\ldots \}$ partition  $A\setminus \{ 0\}$.

There exists a positive integer $n_0$ such that, if $n \geq n_0$ and $n \in hA$, 
then $r_{A,h}(n) \geq \ell$.  Because $A$ is infinite, there exists $a_0 \in A$ 
with $ha_0 \geq n_0$.  Choose $k_0$ such that $a_0 \in A_{k_0}$.

Suppose that $k \geq k_0$ and $A_k \neq \emptyset$.  
Let 
\[
a_k^* = \max(A_k).
\]
Then 
\[
h^{k-1} \leq a_k^* < h^k
\]
and 
\[
A \cap \left[a_k^*+1,h^k \right) = \emptyset.
\]
Consider the integer
\[
ha_k^* \in hA.
\]
Because $a_k^* \geq a_0$, we have  $ha_k^* \geq ha_0 \geq n_0$, 
and so $r_{A,h}(ha_k^*) \geq \ell$.   
It follows that the set $A$ contains a subset 
\[
\{ a_{i,j}:i=1,\ldots, h \text{ and } j = 1,\ldots, \ell -1\}
\]
such that,  for $j=1,2,\ldots, \ell -1$,  
\beq   \label{filaseta:ineq3} 
a_{1,j} < a_{h,j} 
\eeq
\beq   \label{filaseta:ineq4}
a_{1,j} \leq a_{2,j} \leq \cdots  \leq a_{h-1,j}  \leq a_{h,j} 
 \eeq 
\[ 
 a_{1,j} +a_{2,j} + \cdots + a_{h-1,j} + a_{h,j} = ha_k^*
\]
 and
\[
 \left( a_{1,j} , a_{2,j}, \ldots + a_{h-1,j} , a_{h,j} \right) \neq 
  \left( a_{1,j'} , a_{2,j'}, \ldots + a_{h-1,j'} , a_{h,j'} \right) 
\]
for $1 \leq j < j' \leq \ell -1$.  
Moreover,  for $i=1,\ldots, h-1$ and $j = 1, \ldots, \ell -1$, 
the inequality $a_{i,j}  \geq 0$ implies that 
\[
a_{h,j} \leq ha_k^*
\]
for $j =1,2,\ldots, \ell-1$.    
 
Let $b \in A$ and let $J$ be a  subset of $\{1,...,\ell -1\}$ such 
that $a_{h,j} = b$ for all $j \in J$.
If $j \in J$, then 
\[
a_{1,j} +a_{2,j} + \cdots + a_{h-1,j}  = ha_k^* - a_{h,j}   = ha_k^* - b
\]
and so 
 \[
 r_{A,h-1}(ha_k^* - b) \geq |J|.
\] 
Because $A$ is a $B_{h-1,s}$-set, we have
\[
r_{A,h-1}(ha_k^* - b) \leq s
\]
 and so $ |J| \leq s$.
 The pigeonhole principle implies that 
 \[
\left| \left\{  a_{h,j}:j=1,\ldots, \ell -1  \right\}  \right| \geq \frac{\ell - 1}{s}.
 \]
It follows from inequalities~\eqref{filaseta:ineq3} and~\eqref{filaseta:ineq4}  that 
\[
ha_k^*=   a_{1,j} +a_{2,j} + \cdots + a_{h,j}  < ha_{h,j}
\]
and so
\[
a_k^* < a_{h,j}  \leq ha_k^* < h^{k+1}.
\]
Therefore,
\[
a_k^* + 1 \leq a_{h,j}  < h^{k+1}.
\]
Equivalently, 
\[
a_{h,j} \in \left[ a_k^*+1,h^{k+1} \right) =  \left[ a_k^*+1, h^k \right) \cup A_{k+1}
\]
 for $j =1,2,\ldots, \ell-1$.  
Because $A \cap \left[a_k^*+1, h^k \right) = \emptyset$,
we see that 
\[
\{ a_{h,j}  : j=1,2,\ldots, \ell-1 \} \subseteq A_{k+1}
\]
and so 
\[
\left| A_{k+1}\right| \geq \frac{\ell - 1}{s}.
\]
It follows by induction that  $\left| A_k \right| \geq (\ell - 1)/s$ for all $k \geq k_0 + 1$.

Let $x \geq h$, and choose the positive integer $t$ such that 
\[
h^t \leq x < h^{t + 1}.
\]
Because
\[
A \setminus \{ 0\} = \bigcup_{k= 1}^{\infty} A_k
\]
it follows that 
\begin{align*}
A(x) & \geq A\left( h^t \right) \geq  \left(\frac{\ell - 1}{s}\right)(j - k_0 ) \\
& >  \left( \frac{\ell - 1}{s} \right) \left( \frac{\log x}{\log h} - k_0 - 1 \right) 
\end{align*}
for $x \geq h$.  Let  $w_0 =  (\ell - 1) (k_0 + 1)/s$.  
This completes the proof.

\end{proof}

 \def\cprime{$'$} \def\cprime{$'$} \def\cprime{$'$}
\providecommand{\bysame}{\leavevmode\hbox to3em{\hrulefill}\thinspace}
\providecommand{\MR}{\relax\ifhmode\unskip\space\fi MR }
\providecommand{\MRhref}[2]{%
  \href{http://www.ams.org/mathscinet-getitem?mr=#1}{#2}
}
\providecommand{\href}[2]{#2}

\end{document}